\documentclass[11pt]{article}
\usepackage{mathrsfs}
\usepackage{}
\usepackage{amsfonts}
\usepackage{amsmath}
\usepackage{amssymb}

\def\sqr#1#2{{\vcenter{\vbox{\hrule height.#2pt
              \hbox{\vrule width.#2pt height#1pt \kern#1pt \vrule
width.#2pt}
              \hrule height.#2pt}}}}
\def\signed #1{{\unskip\nobreak\hfil\penalty50
              \hskip2em\hbox{}\nobreak\hfil#1
              \parfillskip=0pt \finalhyphendemerits=0 \par}}
\def\endpf{\signed {$\sqr69$}}
\def\3n{\negthinspace \negthinspace \negthinspace }
\def\2n{\negthinspace \negthinspace }
\def\1n{\negthinspace }

\def\no{\noindent}

\def\({\Big (}
\def\){\Big )}
\def\[{\Big[}
\def\]{\Big]}

\def\be{\begin{equation}}
\def\bel{\begin{equation}\label}
\def\ee{\end{equation}}
\def\bea{\begin{eqnarray}}
\def\eea{\end{eqnarray}}
\def\bt{\begin{theorem}}
\def\et{\end{theorem}}
\def\bc{\begin{corollary}}
\def\ec{\end{corollary}}
\def\bl{\begin{lemma}}
\def\el{\end{lemma}}
\def\bp{\begin{proposition}}
\def\ep{\end{proposition}}
\def\br{\begin{remark}}
\def\er{\end{remark}}
\def\ba{\begin{array}}
\def\ea{\end{array}}
\def\bd{\begin{definition}}
\def\ed{\end{definition}}

\newtheorem{lemma}{Lemma}[section]
\newtheorem{remark}{Remark}[section]

\newtheorem{theorem}{Theorem}[section]
\newtheorem{corollary}{Corollary}[section]

\newtheorem{definition}{Definition}[section]
\newtheorem{proposition}{Proposition}[section]

\makeatletter
   
   \@addtoreset{equation}{section}
\makeatother
\allowdisplaybreaks
\begin{document}

\title{\bf Hierarchical control for the semilinear parabolic equations with interior degeneracy }
\author{Hang Gao, Wei Yang  and  Muming Zhang* }

\date{School of Mathematics and Statistics, Northeast Normal
University,\\ Changchun 130024, China.}

\maketitle

\renewcommand{\thefootnote}{\fnsymbol{footnote}}
\footnote[0]{ *Corresponding author (M. Zhang): zhangmm352@nenu.edu.cn.}
\footnote[0]{This work is partially
supported by the NSF of China under grants 12001094, 12001087 and 11971179, and Fundamental Research Funds for the Central Universities under grant 2412020QD027. }

\begin{abstract}
This paper concerns with the hierarchical control of the semilinear parabolic equations with interior degeneracy. By a Stackelberg-Nash strategy, we consider the linear and semilinear system with one leader and two followers. First, for any given leader, we analyze a Nash equilibrium corresponding to a bi-objective optimal control problem. The existence and uniqueness of the Nash equilibrium is proved, and its characterization is given. Then,  we find a leader satisfying the null controllability problem. The key is to establish a new Carleman estimate for a coupled degenerate parabolic system with interior degeneracy.
\end{abstract}

\no{\bf Key Words. }Parabolic equation, Interior degeneracy,
Stackelberg-Nash strategy, Carleman estimate

\no{\bf 2020 Mathematics Subject Classification 2020.} 93B05; 93B07; 35K65; 90C29.

\section{Introduction and main results}

Let $T>0$ and $Q=(0,1)\times(0,T)$. Assume that $\omega,\ \omega_{1},$ and $\omega_{2}$ are three given nonempty open subsets of~$(0,1)$ and $\omega_{i}\cap\omega=\emptyset(i=1,2)$. Denote by  $\chi_{\omega}$ the characteristic function of the set $\omega$. 
We consider  the following degenerate parabolic equation with interior degeneracy:\begin{equation}\label{6.1}
\left\{\begin{array}{ll}
y_{t}-(a(x)y_{x})_{x}+c(x,t)y=G(y)+f\chi_{\omega}+u_{1}\chi_{\omega_{1}}+u_{2}\chi_{\omega_{2}}& \mbox{ in }Q,\\[3mm]
y(0,t)=y(1,t)=0&\mbox{ in }(0,T),\\[3mm]
y(x,0)=y_{0}(x)&\mbox{ in }(0,1),
\end{array}\right.
\end{equation}
where $c\in L^{\infty}(Q),$ $G$ is a locally lipschitz-continuous function,  $y_{0}\in\ L^{2}(0,1)$ is a given initial value, $f\in\ L^{2}(\omega\times(0,T))$~and~$u_{i}\in\ L^{2}(\omega_{i}\times(0,T))(i=1,2)$ are the leader and follower control functions, respectively,  and $y=y(\cdot,\cdot;f,u_{1},u_{2})$ is the state. The function $a$ degenerates at the point $x_{0}\in(0,1)$ with $x_{0}\in\omega$. Assume 
\begin{equation}\label{zmm}
a\in C^{1}([0,1]\backslash\{{x_{0}}\}),\ a>0\ 
\text{in}\ [0,1]
\\\noindent\backslash\{{x_{0}}\},\ a(x_{0})=0\  \text{with}\ x_{0}\in\omega.
\end{equation}

\noindent {\bf Condition 1} Weakly degenerate case:  in addition to (\ref{zmm}), there exists $K\in(0,1)$ such that $(x-x_{0})a'\leq Ka$ in $[0,1]\backslash\{{x_{0}}\}$.

\noindent {\bf Condition 2} Strongly degenerate case: in addition to (\ref{zmm}), $a\in W^{1,\infty}(0,1)$, and there exists ~$K\in[1,2)$ such that $(x-x_{0})a'\leq Ka$~in~$[0,1]\backslash\{{x_{0}}\}$.

For example,  if $a(x)=|x-x_{0}|^{\alpha},$ then $0<\alpha<1$ and $1\leq\alpha<2$ correspond to  weakly degenerate case and strongly degenerate case, respectively.

Degenerate parabolic equations can be used to describe  a wide variety of problems in physics, economics, biology and mathematical finance, for example, boundary layer models, Grushin type models and Fleming-Viot models (see \cite{5,9,20} and the rich references therein). The purpose of  this paper is to investigate 
the hierarchic null controllability of  degenerate parabolic equation (\ref{6.1})  through Stackelberg-Nash strategies. To this aim, 
we define the main cost functional:
\begin{eqnarray}\label{6.2}
J(f)=\frac{1}{2}\int_{\omega\times(0,T)}|f|^{2}dxdt,
\end{eqnarray}
and the secondary cost functionals:
\begin{eqnarray}\label{6.3}
J_{i}(f,u_{1},u_{2})=\frac{\alpha_{i}}{2}\int_{O_{i,d}\times(0,T)}|y-y_{i,d}|^{2}dxdt+\frac{\mu_{i}}{2}\int_{\omega_{i}\times(0,T)}\sigma^{2}|u_{i}|^{2}dxdt, ~i=1,2,
\end{eqnarray}
where $\alpha_{i}~and ~\mu_{i}$~are two positive constants, $y_{i,d}\in L^{2}\left(O_{i,d}\times(0,T)\right)$ are given functions, $\sigma=\sigma(t)$~is a positive function which will be given in (\ref{7.3}), $O_{i,d}\subseteq(0,1)(i=1,2)$~are two observation domains,  $y=y(\cdot,\cdot;f,u_{1},u_{2})$~is the solution of system (\ref{6.1}) corresponding to leader control~$f$~and followers controls~$(u_{1},u_{2})$.
Nash equilibrium pairs and Nash quasi-equilibrium pairs are defined for linear and semilinear cases, respectively.

\begin{definition}\label{d1}
Suppose that~$G\equiv0,$~for any given leader control~$f$, a follower control pair $(\bar{u}_{1},\bar{u}_{2})$ is called a Nash equilibrium pair of~$J_{i}$, if  the following inequalities
\begin{eqnarray}
&J_{1}(f,\bar{u}_{1},\bar{u}_{2})\leq J_{1}(f,u_{1},\bar{u}_{2}),~~\forall u_{1}\in L^{2}(\omega_{1}\times(0,T)), \label{6.4} \\[2mm]
&J_{2}(f,\bar{u}_{1},\bar{u}_{2})\leq J_{2}(f,\bar{u}_{1},u_{2}),~~\forall u_{2}\in L^{2}(\omega_{2}\times(0,T))\label{6.5}
\end{eqnarray}
simultaneously hold.
\end{definition}

\begin{definition}\label{d2}
Suppose that~$G\neq0$,  for any given leader control~$f$,
 if $(\bar{u}_{1},\bar{u}_{2})$ satisfies the following equalities
 \begin{eqnarray}
&\left(\frac{\partial J_{1}}{\partial u_{1}}(f,\bar{u}_{1},\bar{u}_{2}),u_{1}\right)=0,~\forall u_{1}\in L^{2}(\omega_{1}\times(0,T)), \label{6.6}\\[2mm]
&\left(\frac{\partial J_{2}}{\partial u_{2}}(f,\bar{u}_{1},\bar{u}_{2}),u_{2}\right)=0,~\forall u_{2}\in L^{2}(\omega_{2}\times(0,T))\label{6.7}
\end{eqnarray}
where~$\left(\displaystyle\frac{\partial J_{i}}{\partial u_{i}}(f,\bar{u}_{1},\bar{u}_{2}),u_{i}\right)(i=1,2)$~denotes the $G\hat{a}teaux$~differentiation of~$J_{i}$~at~$(f,\bar{u}_{1},\bar{u}_{2})$~along the direction of~$u_{i}$. Then the follower control pair~$(\bar{u}_{1},\bar{u}_{2})$~is said to be a Nash quasi-equilibrium pair of functional~$J_{i}$.
\end{definition}

The main goal of this paper is to study the null controllability  of system~(\ref{6.1})~under Stackelberg-Nash strategies, namely, our two objectives are as follows: 

\noindent 1) For any given leader $f$,  a Nash equilibrium (or quasi-equilibrium) pair of $J_{i}$ exists, denoted as $(\bar{u}_{1}(f),\bar{u}_{2}(f))$.

\noindent 2) There exists a leader control $\bar{f}\in L^{2}(\omega\times(0,T))$ such that \begin{eqnarray}\label{6.8}
J(\bar{f})=\min_{f}J(f,\bar{u}_{1}(f),\bar{u}_{2}(f)),~\forall f\in L^{2}(\omega\times(0,T)),
\end{eqnarray}
and the corresponding solution of  system (\ref{6.1}) satisfies
\begin{eqnarray}\label{6.9}
y(\cdot,T;\bar{f},\bar{u}_{1}(\bar{f}),\bar{u}_{2}(\bar{f}))=0 \ \text{in}\ (0,1).
\end{eqnarray}

Before giving the main result of this paper, we assume that~$a$~satisfies the following condition:
\\\noindent {\bf Condition 3}~The function~$a$~satisfies Condition 1 or Condition 2, and there exists a constant $\rho\in(0,K]$ such that the function $\displaystyle\frac{a(x)}{|x-x_{0}|^{\rho}}$ is nonincreasing in $(0,x_{0})$ and nondecreasing in $(x_{0},1)$, where the constant $K$ appears in Conditions 1 or Condition 2.

Now, the main results in this paper are stated as follows.
\begin{theorem}\label{T1}
Suppose the following conditions hold:
\\\noindent(1)~$O_{1,d}=O_{2,d}$ denoted as $O_{d}$ and $O_{d}\cap\omega\neq\emptyset$, the constants $\mu_{1}$ and $\mu_{2}$ are large enough,
\\\noindent(2)~$x_{0}\in\omega$~and Condition 3 holds,~$G\in W^{1,\infty}(\mathbb{R})$,
\\\noindent(3) $y_{1,d}$ and $y_{2,d}$ satisfy
\begin{eqnarray}\label{6.10}
\int_{O_{d}\times(0,T)}|y_{i,d}|^{2}dxdt<+\infty,~i=1,2.
\end{eqnarray}
Then for any~$y_{0}\in\ L^{2}(0,1)$, there exists a leader control $\bar{f}\in\ L^{2}(\omega\times(0,T))$ and the corresponding Nash equilibrium (or quasi-equilibrium ) pair $(\bar{u}_{1}(\bar{f}),\bar{u}_{2}(\bar{f}))$ such that the solution of system (\ref{6.1}) satisfy (\ref{6.8})-(\ref{6.9}).
\end{theorem}

The hierarchic controls were introduced by Lions in \cite{28,29}, which study the bi-objective control problems for wave equations and heat equations, respectively.  In the past decade, there are a large number of works attributed to the  hierarchic control problem of PDEs. Most of the works dealing with hierarchic control employing the Strackelberg strategies.
In \cite{36}, Stackelberg competition was introduced as a strategy game between several firms. 
 Among the firms, one moves first (called the leader), while the others (named followers) react to the leader' s strategy. When there are many followers, each with a specific optimality objective, Nash equilibrium is the best choice.
  In the past years, there are some works addressing the Stackelberg-Nash controllability problems of linear and semilinear parabolic equations (see, for example, \cite{M4,M2,M3,M6}). 
  We also refer to \cite{nonlinear} for the application of 
 Stackelberg-Nash controllability for $N$-dimensional nonlinear parabolic equations.  
 Recently, The same problem for linear  parabolic equations with dynamic boundary conditions and drift terms was investigated in \cite{arxiv}.   In \cite{coupled}, the Stackelberg-Nash controllability problem of coupled parabolic equations was discussed. However, to our knowledge, very few results are obtainable on the Stackelberg-Nash controllability of degenerate parabolic equations.  In \cite{AAFC}, the authors proved  Stackelberg-Nash controllability for some linear and semilinear  degenerate parabolic equations, whose degenerate point is the origin (the boundary of the domain). Another relevant work concerning the Stackelberg-Nash controllability for nonlinear coupled parabolic equations with boundary degeneracy, refers to \cite{m7}. To the best of our knowledge, Theorem \ref{T1} in this paper is the first result regarding the Stackelberg-Nash controllability  for parabolic equation with interior degeneracy.
 As usual, we first prove the existence and uniqueness of the Nash equilibrium (or quasi-equilibrium). Moreover, its characterization is described. Next, we establish the null controllability of coupled parabolic equations with interior degeneracy by deriving a new Carleman estimate. 
 Controllability and Carleman estimates for degenerate parabolic equations, including one-dimensional and high-dimensional case, have received a lot of attention in the past decade (see \cite{m10,m9,m13,m11,m12,m14,m15,m16} and the rich references therein). 
 In this respect, in \cite{m8}, the author obtained the Carleman estimates for parabolic equations with interior degeneracy in a set of positive measure recently.
We refer to \cite{T} for the Carleman estimates of parabolic equations with interior single point degeneracy. Based on this, we establish a new Carleman estimate for coupled systems.  
 
 The rest of this paper is organized as follows. In Section 2, we present some technical lemmas. In Section 3, we prove the existence and uniqueness of Nash equilibrium and give its characterization. We devote Section 4 to deriving a new Carleman estimate and proving the null controllability result for the linear system. The semilinear case is analyzed in Section 5 by using fixed-point technique.

\medskip

\section{Some technical lemmas}

To prove the main results, some lemmas are provided in this section. First, we introduce some appropriate weighted spaces.

\medskip
(1) Weakly degenerate case: 
\begin{eqnarray*}
H^{1}_{a}(0,1):=\{z\in L^{2}(0,1)|z~is~absolutely~continuous~in~[0,1],~\sqrt{a}z_{x}\in L^{2}(0,1)~and~z(0)=z(1)=0\}
\end{eqnarray*}
and$$H^{2}_{a}(0,1):=\{z\in H^{1}_{a}(0,1)|az_{x}\in H^{1}(0,1)\}.$$

(2) Strongly degenerate case: 
\begin{eqnarray*}
\begin{array}{rl}
&H^{1}_{a}(0,1):=\Big{\{}z\in L^{2}(0,1)|z~is~locally~absolutely~continuous~in~[0,x_{0})\cup(x_{0},1],~\sqrt{a}z_ {x}\in L^{2}(0,1)\\[5mm]&\quad\quad\quad\quad\quad~and~z(0)=z(1)=0\Big{\}}
\end{array}
\end{eqnarray*}
and$$H^{2}_{a}(0,1):=\{z\in H^{1}_{a}(0,1)|az_{x}\in H^{1}(0,1)\}.$$
The associated norms for two cases whose are
$$\|z\|^{2}_{H^{1}_{a}(0,1)}:=\|z\|^{2}_{L^{2}(0,1)}+\|\sqrt{a}z_{x}\|^{2}_{L^{2}(0,1)}$$
and
$$\|z\|^{2}_{H^{2}_{a}(0,1)}:=\|z\|^{2}_{H^{1}_{a}(0,1)}+\|(az_{x})_{x}\|^{2}_{L^{2}(0,1)},$$
respectively. 
The conjugate space of~$H^{1}_{a}(0,1)$ is denoted by $H^{*}_{a}(0,1)$, and 
$$\|z^{*}\|_{H^{*}_{a}(0,1)}:=\sup_{\|z\|_{H^{1}_{a}(0,1)}=1}\langle z,z^{*}\rangle_{H^{1}_{a}(0,1),H^{*}_{a}(0,1)}.$$
We define the operator $Az:=(az_{x})_{x}$  for all $z\in D(A)$, where  $D(A)=H^{2}_{a}(0,1)$.

\medskip

Next, we recall a konwn  Carleman estimate (see \cite[Lemma 4.1]{T}) and a Hardy-Poincar$\acute{e}$ inequality (see \cite[Proposition 2.3]{T}) for
  parabolic operators with interior degeneracy. We introduce an auxiliary function:
\begin{eqnarray}\label{7.1}
\phi(x,t)=\theta(t)\psi(x),
\end{eqnarray}
where
\begin{eqnarray}\label{7.2}
\begin{array}{rl}
&\theta(t)=\displaystyle\frac{1}{[t(T-t)]^{4}},~\psi(x)=c_{1}\left[\displaystyle\int_{x_{0}}^{x}\frac{y-x_{0}}{a(y)}dy-c_{2}\right],\\[5mm]&c_{2}>\max\left[\displaystyle\frac{(1-x_{0})^{2}}{a(1)(2-K)},~\displaystyle\frac{x^{2}_{0}}{a(0)(2-K)}\right],~c_{1}>0, \end{array}
\end{eqnarray}
and $K$~is the constant that appears in Condition~1~or Condition~2. It is easy to check that $\psi(x)<0$ for all $x\in[0,1]$ and $\psi\geq-c_{1}c_{2}$~(See~(3.4)~in~[20]).  Set
\begin{eqnarray}\label{7.3}
\phi^{*}(t)=\min_{x\in[0,1]}\phi(x,t),\ \sigma(t)=e^{\frac{-s\bar{\phi}^{*}}{2}}.
\end{eqnarray}

\begin{lemma}\label{l1}
Assume Condition 3 holds and~$x_{0}\in\omega$. Then there exist two positive constants~$C$ and~$s_{0}$ such that  for any $s\geq s_{0}$  and any solution $\vartheta$ 
of 
\begin{equation}\label{7.4}
\left\{\begin{array}{ll}
\vartheta_{t}+(a\vartheta_{x})_{x}=g& \mbox{ in }Q,\\[3mm]
\vartheta(0,t)=\vartheta(1,t)=0&\mbox{ in }(0,T),\\[3mm]
\vartheta(x,T)=\vartheta_{T}(x)&\mbox{ in }(0,1),
\end{array}\right.
\end{equation}
it holds that 
$$\int_{Q}\left[s\theta a\vartheta_{x}^{2}+s^{3}\theta^{3}\frac{(x-x_{0})^{2}}{a}\vartheta^{2}\right]e^{2s\phi}dxdt\leq C\left(\int_{\omega\times(0,T)}s^{3}\theta^{3}e^{2s\phi}\vartheta^{2}dxdt+\int_{Q}e^{2s\phi}g^{2}dxdt\right),$$
where $\vartheta_{T}\in D(A^{2}):=\{z\in D(A)|Az\in D(A)\}$,  $\phi$ and $\theta$ are given in (\ref{7.1}) and (\ref{7.2}).
\end{lemma}

\begin{remark}
The Carleman estimate for homogeneous  degenerate parabolic equations is presented in \cite[Lemma 4.1]{T}.  By improving the proof process of \cite[Lemma 4.1]{T},  Lemma \ref{l1} can be obtained directly for non-homogeneous  degenerate parabolic equations. 
\end{remark}

\begin{lemma}\label{l2}
~Assume that~$p\in C([0,1]),p>0$ in~$[0,1]\backslash\{x_{0}\}$, $p(x_{0})=0$, and there exists~$q\in (1,2)$ such that the function~$\displaystyle\frac{p(x)}{|x-x_{0}|^{q}}$ is nonincreasing in $(0,x_{0})$, and nondecreasing in $(x_{0},1)$. Then, there exists a constant~$C>0$~such that for any function~$w$, which is locally
absolutely continuous on~$[0,x_{0})\cup(x_{0},1]$~and satisfies
$$w(0)=w(1)=0 \ \text{and}\ \int_{0}^{1}p(x)|w_{x}(x)|^{2}dx<+\infty,$$
it holds that
\begin{eqnarray}\label{7.5}
\int_{0}^{1}\frac{p(x)}{(x-x_{0})^{2}}w^{2}(x)dx\leq C\int_{0}^{1}p(x)|w_{x}(x)|^{2}dx.
\end{eqnarray}
\end{lemma}

\section{Nash equilibrium}

In this section, we consider the linear case, i.e., $G(y)\equiv0$ in equation (\ref{6.1}). 
The existence and uniqueness of  Nash equilibrium pair will be proved first, then its explicit expression will be given.
In the remainder of this paper, we specify that $i$ includes both cases $i=1$ and $i=2$, 
and $C$ denotes a positive constant depending only on ~$(0,1),O_{d},\omega,T,\alpha_{i},\mu_{i}(i=1,2),$ and $\|c\|_{L^{\infty}(Q)}$,  which may be different from one line to another.

\subsection{Existence and uniqueness of the Nash equilibrium pair}

In this subsection, we will prove the existence and uniqueness of Nash equilibrium pairs, for any given leader control.

 Define the function space
$$S_{i}=L^{2}(\omega_{i}\times(0,T)), ~i=1,2,\quad S:=S_{1}\times S_{2},$$
  and introduce the operator $\Lambda_{i}: S_{i} \rightarrow L^{2}(Q))$,
 by $\Lambda_{i}u_{i}=y^{i}$, where~$y^{i}$~is the solution of the following system
\begin{equation}\label{8.1}
\left\{\begin{array}{ll}
y^{i}_{t}-(ay^{i}_{x})_{x}+cy^{i}=u_{i}\chi_{\omega_{i}}& \mbox{ in }Q,\\[3mm]
y^{i}(0,t)=y^{i}(1,t)=0&\mbox{ in }(0,T),\\[3mm]
y^{i}(x,0)=0&\mbox{ in }(0,1).
\end{array}\right.
\end{equation}
Then, for any~$f\in L^{2}(\omega\times(0,T))$, we can write the solution of (\ref{6.1}) as follows
$$y=y^{1}+y^{2}+\zeta(f)=\Lambda_{1}u_{1}+\Lambda_{2}u_{2}+\zeta(f),$$
where $\zeta(f)$ is the solution of the following system
\begin{equation}\label{8.2}
\left\{\begin{array}{ll}
\zeta_{t}-(a\zeta_{x})_{x}+c\zeta=f\chi_{\omega}& \mbox{ in }Q,\\[3mm]
\zeta(0, t)=\zeta(1,t)=0&\mbox{ in }(0,T),\\[3mm]
\zeta(x,0)=y_{0}(x)&\mbox{ in }(0,1).
\end{array}\right.
\end{equation}
Therefore, the functional (\ref{6.3}) can be expressed as follows
$$J_{i}(f,u_{1},u_{2})=\frac{\alpha_{i}}{2}\int_{O_{i,d}\times(0,T)}|\Lambda_{1}u_{1}+\Lambda_{2}u_{2}-\tilde{y}_{i,d}|^{2}dxdt+\frac{\mu_{i}}{2}\int_{\omega_{i}\times(0,T)}\sigma^{2}|u_{i}|^{2}dxdt,$$
where~$\tilde{y}_{i,d}=y_{i,d}-\zeta(f)|_{O_{i,d}}.$

\medskip
Notice that $J_{1}$ and $J_{2}$ are strictly convex functionals,  by (\ref{6.4})-(\ref{6.5}), we  obtain that $(\bar{u}_{1},\bar{u}_{2})$~is a Nash equilibrium pair of $(J_{1},J_{2})$  if and only if 
\begin{eqnarray}\label{8.3}
\left(\frac{\partial J_{1}}{\partial u_{1}}(f,\bar{u}_{1},\bar{u}_{2}),u_{1}\right)=0,~\forall u_{1}\in S_{1}\ \text{and}\ 
\left(\frac{\partial J_{2}}{\partial u_{2}}(f,\bar{u}_{1},\bar{u}_{2}),u_{2}\right)=0,~\forall u_{2}\in S_{2},
\end{eqnarray}
where~$\left(\displaystyle\frac{\partial J_{1}}{\partial u_{1}}(f,\bar{u}_{1},\bar{u}_{2}),u_{1}\right)$~denotes the G$\hat{a}$teaux~differentiation of~$J_{1}$~at~$(f,\bar{u}_{1},\bar{u}_{2})$~along the direction of~$u_{1}$,  i.e.,
$(\bar{u}_{1},\bar{u}_{2})$ is a Nash equilibrium pair if and only if for any~$(u_{1},u_{2})\in S$, it holds that
\begin{eqnarray}\label{8.4}
\mu_{i}\int_{\omega_{i}\times(0,T)}\sigma^{2}\bar{u}_{i}u_{i}dxdt+\alpha_{i}\int_{O_{i,d}\times(0,T)}(\Lambda_{1}\bar{u}_{1}+\Lambda_{2}\bar{u}_{2}-\tilde{y}_{i,d})\cdot\Lambda_{i}u_{i}dxdt=0,  ~i=1,2.
\end{eqnarray}
That is to say, for any ~$(u_{1},u_{2})\in S$, we get
\begin{eqnarray*}
\mu_{i}(\sigma^{2}\bar{u}_{i},u_{i})_{S_{i}}+\alpha_{i}\left(\Lambda_{i}^{*}\left[(\Lambda_{1}\bar{u}_{1}+\Lambda_{2}\bar{u}_{2})|_{O_{i,d}}-\tilde{y}_{i,d}\right],u_{i}\right)_{S_{i}}=0, ~i=1,2,
\end{eqnarray*}
where~$(\cdot,\cdot)_{S_{i}}$~denotes the inner product of $S_{i}$, and~$\Lambda_{i}^{*}\in \mathcal{L}\left(L^{2}(Q),S_{i}\right)$~is the adjoint operator of $\Lambda_{i}$. 
This implies 
$$\mu_{i}\sigma^{2}\bar{u}_{i}+\alpha_{i}\Lambda_{i}^{*}\left[(\Lambda_{1}\bar{u}_{1}+\Lambda_{2}\bar{u}_{2})|_{O_{i,d}}\right]=\alpha_{i}\Lambda_{i}^{*}\tilde{y}_{i,d}, ~i=1,2.$$

Now, for any~$u=(u_{1},u_{2})$, we define the operator~$H=(H_{1},H_{2})$~as follows, where~$H_{1}\in\mathcal{L}(S,S_{1}),~H_{2}\in\mathcal{L}(S,S_{2})$,
$$H_{i}u=\mu_{i}\sigma^{2}u_{i}+\alpha_{i}\Lambda_{i}^{*}\left[(\Lambda_{1}u_{1}+\Lambda_{2}u_{2})|_{O_{i,d}}\right],~i=1,2.$$
Therefore, $\bar{u}=(\bar{u}_{1},\bar{u}_{2})$~is a Nash equilibrium pair equivalent to
\begin{eqnarray}\label{8.5}
H\bar{u}=(H_{1}\bar{u},H_{2}\bar{u})=(\alpha_{1}\Lambda_{1}^{*}\tilde{y}_{1,d},~\alpha_{2}\Lambda_{2}^{*}\tilde{y}_{2,d}).
\end{eqnarray}

By Cauchy inequality with~$\varepsilon$, we have
\begin{eqnarray}\label{8.6}
\begin{array}{rl}
&(Hu,u)_{S}=(H_{1}u,u_{1})_{S_{1}}+(H_{2}u,u_{2})_{S_{2}}\\[5mm]&\quad\quad\quad~~=\displaystyle\sum_{i=1}^{2}\mu_{i}\|\sigma u_{i}\|^{2}_{L^{2}(\omega_{i}\times(0,T))}+\alpha_{1}(\Lambda_{1}u_{1}+\Lambda_{2}u_{2},\Lambda_{1}u_{1})_{L^{2}(O_{1,d}\times(0,T))}\\[5mm]&\quad~\quad\quad\quad~+\alpha_{2}(\Lambda_{1}u_{1}+\Lambda_{2}u_{2},\Lambda_{2}u_{2})_{L^{2}(O_{2,d}\times(0,T))}.\\[2mm]
&\quad\quad\quad\quad\geq \mu_{1}\|\sigma_{0}^{2}u_{1}\|^{2}_{S_{1}}\!\!+\!\mu_{2}\|\sigma_{0}^{2}u_{2}\|^{2}_{S_{2}}\!\!+\!\alpha_{1}\|\Lambda_{1}u_{1}\|^{2}_{L^{2}(O_{1,d}\times(0,T))}\!\!-\!\displaystyle\frac{\varepsilon}{2}\alpha_{1}\|\Lambda_{2}u_{2}\|_{L^{2}(O_{1,d}\times(0,T))}\\[5mm]&\quad\quad\quad\quad~-\!\displaystyle\frac{1}{2\varepsilon}\alpha_{1}\|\Lambda_{1}u_{1}\|_{L^{2}(O_{1,d}\times(0,T))}\!\!+\!\alpha_{2}\|\Lambda_{2}u_{2}\|^{2}_{L^{2}(O_{2,d}\times(0,T))}\!\!-\!\displaystyle\frac{\varepsilon}{2}\alpha_{2}\|\Lambda_{1}u_{1}\|_{L^{2}(O_{2,d}\times(0,T))}\\[5mm]
&\quad\quad\quad\quad~-\displaystyle\frac{1}{2\varepsilon}\alpha_{2}\|\Lambda_{2}u_{2}\|_{L^{2}(O_{2,d}\times(0,T))}\\[5mm]&\quad\quad~~~~\geq \mu_{1}\|\sigma_{0}^{2}u_{1}\|^{2}_{S_{1}}\!+\!\mu_{2}\|\sigma_{0}^{2}u_{2}\|^{2}_{S_{2}}\!-\!\displaystyle\frac{\alpha_{1}}{4}\|\Lambda_{2}\chi_{O_{1,d}}\|^{2}_{S_{1,d}}\|u_{2}\|^{2}_{S_{2}}\!-\!\displaystyle\frac{\alpha_{2}}{4}\|\Lambda_{1}\chi_{O_{2,d}}\|^{2}_{S_{2,d}}\|u_{1}\|^{2}_{S_{1}},
\end{array}
\end{eqnarray}
where~$\sigma_{0}=\displaystyle\min_{t\in[0,T]}\sigma(t)$, $\varepsilon=1/2$,~$\|\cdot\|_{S_{i,d}}(i=1,2)$~denotes the norm in the space~$\mathcal{L}(S_{3-i},L^{2}(O_{i,d}\times(0,T)))$. Then, for~$\mu_{1},\mu_{2}$~sufficiently large satisfying
$$4\mu_{1}\sigma_{0}^{2}>\alpha_{2}\|\Lambda_{1}\chi_{O_{2,d}}\|^{2}_{S_{2,d}}, ~4\mu_{2}\sigma_{0}^{2}>\alpha_{1}\|\Lambda_{2}\chi_{O_{1,d}}\|^{2}_{S_{1,d}}$$
we obtain
\begin{eqnarray}\label{8.7}
(Hu,u)_{S}\geq\beta\|u\|^{2}_{S}, ~\beta=\min_{i=1,2}\left\{\mu_{i}\sigma_{0}^{2}-\frac{\alpha_{3-i}}{4}\|\Lambda_{i}\chi_{O_{3-i.d}}\|^{2}_{S_{3-i,d}}\right\}>0.
\end{eqnarray}

On the other hand, we define the functional $b(u,v):S\times S\longrightarrow \mathbb{R}$
by $b(u,v):=(Hu,v)_{S}.$
By the definition of~$H$ and (\ref{8.7}), we see that~$b$~is a coercive bounded bilinear form. 
Take~$r=(\alpha_{1}\Lambda_{1}^{*}\tilde{y}_{1,d},\alpha_{2}\Lambda_{2}^{*}\tilde{y}_{2,d})$,~define the functional~$\hat{F}(v)=(r,v)_{S}$, by Lax-Milgram theorem, we obtain that there exists  a unique $\bar{u}\in S$ such that
$b(\bar{u},v)=(U,v)_{S},~\forall v\in S$. Moreover,
\begin{eqnarray}\label{8.8}
\begin{array}{rl}
&\|\bar{u}\|_{S}\leq C\|r\|_{S}\leq C\left(\|\tilde{y}_{1,d}\|_{L^{2}(O_{1,d}\times(0,T))}+\|\tilde{y}_{2,d}\|_{L^{2}(O_{2,d}\times(0,T))}\right)\leq C\|\zeta(f)\|_{L^{2}(Q)}\\[5mm]&\quad\quad~\leq C\left(\|y_{0}\|_{L^{2}(0,1)}+\|f\|_{_{L^{2}(\omega\times(0,T))}}\right).
\end{array}
\end{eqnarray}
Thus, $H\bar{u}=r$, this means that (\ref{8.5}) holds.  The existence and uniqueness of Nash equilibrium pairs of $(J_{1},J_{2})$ is proved.

\subsection{Characterization of the Nash equilibrium pair}

This subsection will provide the explicit expression for a Nash equilibrium pair, which is stated as follows.

\begin{proposition}\label{p1}
Assume that $\mu_{i}(i=1,2)$ is sufficiently large, then for any leader control $f\in\ L^{2}(\omega\times(0,T))$, there exists a unique Nash equilibrium pair $(\bar{u}_{1},\bar{u}_{2})\in S$ of $(J_{1},J_{2})$ such that (\ref{6.4}) and (\ref{6.5}) hold. Moreover,
\begin{eqnarray}\label{8.10}
\bar{u}_{i}=-\frac{1}{\mu_{i}}\sigma^{-2}q^{i}\chi_{\omega_{i}},\quad i=1,2,
\end{eqnarray}
where~$(y,q^{1},q^{2})$~is the solution of the following coupled systems corresponding to the leader control $f$:
\begin{equation}\label{8.11}
\left\{\begin{array}{ll}
y_{t}-(ay_{x})_{x}+cy=f\chi_{\omega}-\frac{1}{\mu_{1}}\sigma^{-2}q^{1}\chi_{\omega_{1}}-\frac{1}{\mu_{2}}\sigma^{-2}q^{2}\chi_{\omega_{2}}& \mbox{ in }Q,\\[3mm]
-q^{i}_{t}-(aq^{i}_{x})_{x}+cq^{i}=\alpha_{i}(y-y_{i,d})\chi_{O_{i,d}}& \mbox{ in }Q,\\[3mm]
y(0,t)=y(1,t)=0,~q^{i}(0,t)=q^{i}(1,t)=0& \mbox{ in }(0,T),\\[3mm]
y(x,0)=y_{0}(x),~q^{i}(x,T)=0& \mbox{ in }(0,1).
\end{array}\right.
\end{equation}
\end{proposition}
\noindent{\bf Proof.}
Notice that $y=\Lambda_{1}\bar{u}_{1}+\Lambda_{2}\bar{u}_{2}+\zeta(f), $ and $\tilde{y}_{i,d}=y_{i,d}-\zeta(f)|_{O_{i,d}}$, by (\ref{8.4}), we know that $(\bar{u}_{1},\bar{u}_{2})$~is a Nash equilibrium pair of$(J_{1},J_{2})$ if and only if
\begin{eqnarray}\label{8.9}
\mu_{i}\displaystyle\int_{\omega_{i}\times(0,T)}\sigma^{2}\bar{u}_{i}u_{i}dxdt+\alpha_{i}\displaystyle\int_{O_{i,d}\times(0,T)}(y-y_{i,d})\cdot y^{i}dxdt=0,~\forall u_{i}\in S_{i},~i=1,2,
\end{eqnarray}
where~$y^{i}$~ is the solution of system (\ref{8.1}). 

On the other hand, we consider the following adjoint system of (\ref{8.1}):
\begin{equation}\label{8.12}
\left\{\begin{array}{ll}
-q^{i}_{t}-(aq^{i}_{x})_{x}+cq^{i}=\alpha_{i}(y-y_{i,d})\chi_{O_{i,d}}& \mbox{ in }Q,\\[3mm]
q^{i}(0,t)=q^{i}(1,t)=0&\mbox{ in }(0,T),\\[3mm]
q^{i}(x,T)=0&\mbox{ in }(0,1),
\end{array}\right.
\end{equation}
where~$y$~is the solution of system (\ref{6.1}),~$y_{i,d}$~are the given objective functions. Multiplying both sides of the first equation of (\ref{8.12}) by $y^{i}$ and integrating it on $Q$, we get
$$\int_{Q}q^{i}(y^{i}_{t}-(ay^{i}_{x})_{x}+cy^{i})dxdt=\int_{Q}\alpha_{i}(y-y_{i,d})\chi_{O_{i,d}}y^{i}dxdt.$$
Combining this with (\ref{8.1}) and~(\ref{8.9}), we have
$$\int_{\omega_{i}\times(0,T)}q^{i}u_{i}dxdt=-\mu_{i}\int_{\omega_{i}\times(0,T)}\sigma^{2}\bar{u}_{i}u_{i}dxdt,~\forall u_{i}\in S_{i},~i=1,2.$$
This implies
$$q^{i}\chi_{\omega_{i}}+\mu_{i}\sigma^{2}\bar{u}_{i}=0,$$
which completes the proof.\endpf

\section{Null Controllability and Carleman estimate of the linear system}

After proving that there exists a Nash equilibrium for each leader control $f$, we only need to find a control ~$\bar{f}\in L^{2}(\omega\times(0,T))$ for system (\ref{8.11}) such that  \begin{eqnarray}\label{8.14}
J(\bar{f})=\min_{f}J(f),~\forall f\in L^{2}(\omega\times(0,T))
\end{eqnarray}
subject to the following null controllability condition
\begin{eqnarray}\label{8.13}
y(\cdot,T;\bar{f})=0 \mbox{~in~}(0,1),
\end{eqnarray}
where $y$ is the solution of (\ref{8.11}).

In order to accomplish this,  we consider the following adjoint of (\ref{8.11}):
\begin{equation}\label{8.15}
\left\{\begin{array}{ll}
-\varphi_{t}-(a\varphi_{x})_{x}+c\varphi=\alpha_{1}\psi^{1}\chi_{O_{1,d}}+\alpha_{2}\psi^{2}\chi_{O_{2,d}}& \mbox{ in }Q,\\[3mm]
\psi^{i}_{t}-(a\psi^{i}_{x})_{x}+c\psi^{i}=-\frac{1}{\mu_{i}}\sigma^{-2}\varphi\chi_{\omega_{i}}& \mbox{ in }Q,\\[3mm]
\varphi(0,t)=\varphi(1,t)=0,~\psi^{i}(0,t)=\psi^{i}(1,t)=0& \mbox{ in }(0,T),\\[3mm]
\varphi(x,T)=\varphi_{T},~\psi^{i}(x,0)=0& \mbox{ in }(0,1),
\end{array}\right.
\end{equation}
where~$\varphi_{T}\in L^{2}(0,1)$~is the given initial value. 
Assume that $O_{1,d}=O_{2,d}=O_{d}$, (\ref{8.15}) can be simplified as follows
\begin{equation}\label{8.16}
\left\{\begin{array}{ll}
-\varphi_{t}-(a\varphi_{x})_{x}+c\varphi=(\alpha_{1}\psi^{1}+\alpha_{2}\psi^{2})\chi_{O_{d}}& \mbox{ in }Q,\\[3mm]
\psi^{i}_{t}-(a\psi^{i}_{x})_{x}+c\psi^{i}=-\frac{1}{\mu_{i}}\sigma^{-2}\varphi\chi_{\omega_{i}}& \mbox{ in }Q,\\[3mm]
\varphi(0,t)=\varphi(1,t)=0,~\psi^{i}(0,t)=\psi^{i}(1,t)=0& \mbox{ in }(0,T),\\[3mm]
\varphi(x,T)=\varphi_{T},~\psi^{i}(x,0)=0& \mbox{ in }(0,1).
\end{array}\right.
\end{equation}
We transform the null controllability problem of (\ref{8.11}) into a suitable observability problem for (\ref{8.16}).

\begin{proposition}\label{ip2}
If there exists a constant~$C>0$~such that the  corresponding solution~$(\varphi,\psi_{i})$~of (\ref{8.16}) satisfies
\begin{eqnarray}\label{8.17}
\displaystyle\int_{0}^{1}|\varphi(x,0)|^{2}dx+\displaystyle\sum_{i=1}^{2}\int_{Q}|\psi^{i}|^{2}dxdt\leq C\int_{\omega\times(0,T)}|\varphi|^{2}dxdt
\end{eqnarray}
 for any~$\varphi_{T}\in L^{2}(0,1)$, then system  (\ref{8.11}) is null controllable.
\end{proposition}

\no{\bf Proof.} Define the following functional for any~$\varphi_{T}\in L^{2}(0,1)$:
\begin{eqnarray}\label{8.18}
F(\varphi_{T})=\frac{1}{2}\int_{\omega\times(0,T)}|\varphi|^{2}dxdt+\int_{0}^{1}y_{0}(x)\varphi(x,0)dx-\displaystyle\sum_{i=1}^{2}\int_{O_{i,d}\times(0,T)}\alpha_{i}\psi^{i}y_{i,d}dxdt,
\end{eqnarray}
where~$(\varphi,\psi^{i})$ is the solution of (\ref{8.16}).
It is clear that the functional $F$ is continuous and strictly convex. In addition, by Cauchy inequality with~$\varepsilon$~and (\ref{8.17}), we get
\begin{eqnarray*}
\begin{array}{rl}
F(\varphi_{T})\geq\displaystyle\frac{1-\varepsilon C}{2}\int_{\omega\times(0,T)}\varphi^{2}dxdt-\frac{1}{2\varepsilon}\left[\int_{0}^{1}|y_{0}|^{2}dx+\displaystyle\sum_{i=1}^{2}\int_{O_{i,d}\times(0,T)}\alpha_{i}^{2}|y_{i,d}|^{2}dxdt\right].
\end{array}
\end{eqnarray*}
Take~$\varepsilon=\displaystyle\frac{1}{2C}$, we have
\begin{eqnarray*}
\begin{array}{rl}
F(\varphi_{T})\geq\displaystyle\frac{1}{4}\int_{\omega\times(0,T)}\varphi^{2}dxdt-C\left[\int_{0}^{1}|y_{0}|^{2}dx+\displaystyle\sum_{i=1}^{2}\int_{O_{i,d}\times(0,T)}\alpha_{i}^{2}|y_{i,d}|^{2}dxdt\right],
\end{array}
\end{eqnarray*}
combining with~(\ref{6.10}), we find that $F(\varphi_{T})$~is coercive. Therefore,  there exists an extremal function $\bar{\varphi}_{T}$ of $F$ satisfying
~$F(\bar{\varphi}_{T})=\displaystyle\min_{\varphi_{T}\in L^{2}(0,1)} F(\varphi_{T}),$
it follows that 
\begin{eqnarray}\label{8.19}
\begin{array}{rl}
\displaystyle\int_{\omega\times(0,T)}\bar{\varphi}\varphi dxdt+\displaystyle\int_{0}^{1}y_{0}(x)\varphi(x,0)dx-\displaystyle\sum_{i=1}^{2}\displaystyle\int_{O_{i,d}\times(0,T)}\alpha_{i}\psi^{i}y_{i,d}dxdt=0
\end{array}
\end{eqnarray}
for any~$\varphi_{T}\in L^{2}(0,1)$,
where~$(\bar{\varphi}, \bar{\psi}^{i}), (\varphi,\psi^{i})$ are the solutions of~(\ref{8.16}) corresponding to $\bar{\varphi}_{T}, \varphi_{T}$, respectively. 

On the other hand, multiplying both sides of the first two equations of  (\ref{8.11})  by $\varphi$ and $\psi^i(i=1,2)$, respectively,  and integrating them on $Q$, by (\ref{8.11}) and (\ref{8.15}), we conclude that
$$\int_{0}^{1}y(x,T)\varphi_{T}dx-\int_{0}^{1}y_{0}(x)\varphi(x,0)dx=\int_{\omega\times(0,T)}f\varphi dxdt-\displaystyle\sum_{i=1}^{2}\int_{O_{i,d}\times(0,T)}\alpha_{i}\psi^{i}y_{i,d}dxdt.$$
This implies that, if  (\ref{8.13}) holds, then
\begin{eqnarray}\label{p4.12}
\int_{\omega\times(0,T)}f\varphi dxdt-\displaystyle\sum_{i=1}^{2}\int_{O_{i,d}\times(0,T)}\alpha_{i}\psi^{i}y_{i,d}dxdt+\int_{0}^{1}y_{0}(x)\varphi(x,0)dx=0.
\end{eqnarray}
Combining (\ref{8.19}) with (\ref{p4.12}), we obtain that (\ref{8.13}) holds if we take 
 the control $\bar{f}=\bar{\varphi}\chi_{\omega}$. Moreover, by  the $L^2$-estimate for the parabolic equations (\ref{8.16}), we get that 
\begin{eqnarray}\label{8.20}
\|\bar{f}\|_{L^{2}(\omega\times(0,T))}=\|\bar{\varphi}\chi_{\omega}\|_{L^{2}(\omega\times(0,T))}\leq C\|\bar{\varphi}_{T}\|_{L^{2}(0,1)}
\end{eqnarray}
for sufficiently small positive constants  $\mu_{i}(i=1,2)$,
the proof is completed.
\endpf

\medskip

Next, we only need to prove the observability estimate (\ref{8.17}) for (\ref{8.16}). We set $\gamma=\alpha_{1}\psi^{1}+\alpha_{2}\psi^{2}$, then (\ref{8.16}) can be simplified as 
\begin{equation}\label{8.21}
\left\{\begin{array}{ll}
-\varphi_{t}-(a\varphi_{x})_{x}+c\varphi=\gamma\chi_{O_{d}}& \mbox{ in }Q,\\[3mm]
\gamma_{t}-(a\gamma_{x})_{x}+c\gamma=-\left(\frac{\alpha_{1}}{\mu_{1}}\chi_{\omega_{1}}+\frac{\alpha_{2}}{\mu_{2}}\chi_{\omega_{2}}\right)\sigma^{-2}\varphi& \mbox{ in }Q,\\[3mm]
\varphi(0,t)=\varphi(1,t)=0,~\gamma(0,t)=\gamma(1,t)=0& \mbox{ in }(0,T),\\[3mm]
\varphi(x,T)=\varphi_{T},~\gamma(x,0)=0& \mbox{ in }(0,1).
\end{array}\right.
\end{equation}
In order to prove (\ref{8.17}), we derive the following global Carleman estimate for system (\ref{8.21}).

\begin{proposition}\label{p3}
Assume that $O_{d}\cap\omega\neq\emptyset,~x_{0}\in\omega$~and Condition 3 holds. Then one can find two positive constants~$C$ and~$s_{0}$ such that for any $s\geq s_{0}$ and any solution $(\varphi, \gamma)$ of (\ref{8.21}), it holds that
\begin{eqnarray}\label{8.22}
I(\varphi)+I(\gamma)\leq C\displaystyle\int_{\omega\times(0,T)}s^{7}\theta^{7}e^{2s\phi}\varphi^{2}dxdt,
\end{eqnarray}
where $I(v):=\displaystyle\int_{Q}\left[s\theta av_{x}^{2}+s^{3}\theta^{3}\displaystyle\frac{(x-x_{0})^{2}}{a}v^{2}\right]e^{2s\phi}dxdt$,~$\phi$~and~$\theta$ are the functions  given in~(\ref{7.1})~and~(\ref{7.2}).
\end{proposition}

\no{\bf Proof.}
We choose $\omega^{\prime}$ to be a nonempty open subset of $(0,1)$ such that $\overline{\omega^{\prime}}\subset O_{d}\cap\omega$. Applying Lemma~\ref{7.1}~to the first two equations in (\ref{8.21}), respectively, we obtain that
\begin{eqnarray}\label{8.23}
\begin{array}{rl}
&I(\varphi)+I(\gamma)\leq C\displaystyle\int_{\omega^{\prime}\times(0,T)}s^{3}\theta^{3}e^{2s\phi}\varphi^{2}dxdt+C\displaystyle\int_{\omega^{\prime}\times(0,T)}s^{3}\theta^{3}e^{2s\phi}\gamma^{2}dxdt\\[5mm]&\quad\quad\quad\quad\quad~+C\displaystyle\int_{Q}e^{2s\phi}|\gamma\chi_{O_{d}}|^{2}dxdt+C\displaystyle\int_{Q}e^{2s\phi}|c\varphi|^{2}dxdt\\[5mm]&\quad\quad\quad\quad\quad~+C\displaystyle\int_{Q}e^{2s\phi}\left|-\left(\frac{\alpha_{1}}{\mu_{1}}\chi_{\omega_{1}}+\frac{\alpha_{2}}{\mu_{2}}\chi_{\omega_{2}}\right)\sigma^{-2}\varphi\right|^{2}dxdt+C\displaystyle\int_{Q}e^{2s\phi}|c\gamma|^{2}dxdt.
\end{array}
\end{eqnarray}
By Young inequality, it holds that
\begin{eqnarray}\label{8.24}
\begin{array}{rl}
&\displaystyle\int_{Q}e^{2s\phi}|\gamma\chi_{O_{d}}|^{2}dxdt=\displaystyle\int_{Q}\left(e^{2s\phi}\frac{a^{\frac{1}{3}}}{|x-x_{0}|^{\frac{2}{3}}}|\gamma\chi_{O_{d}}|^{2}\right)^{\frac{3}{4}}\left(e^{2s\phi}\frac{|x-x_{0}|^{2}}{a}|\gamma\chi_{O_{d}}|^{2}\right)^{\frac{1}{4}}dxdt\\[5mm]&\quad~~\quad\quad\quad\quad\quad\quad\quad\leq\displaystyle\frac{3}{4}\displaystyle\int_{Q}e^{2s\phi}\frac{a^{\frac{1}{3}}}{|x-x_{0}|^{\frac{2}{3}}}\gamma^{2}dxdt+\displaystyle\frac{1}{4}\displaystyle\int_{Q}e^{2s\phi}\frac{|x-x_{0}|^{2}}{a}\gamma^{2}dxdt.
\end{array}
\end{eqnarray}
Take~$p(x)=a(x)\left(\displaystyle\frac{|x-x_{0}|^{2}}{a}\right)^{\frac{2}{3}},$~then $\displaystyle\frac{p(x)}{(x-x_{0})^{2}}=\displaystyle\frac{a^{\frac{1}{3}}}{|x-x_{0}|^{\frac{2}{3}}}$. Combining Condition~3~and Lemma~\ref{l2}, indicates
\begin{eqnarray}\label{8.25}
\begin{array}{rl}
&\displaystyle\int_{Q}e^{2s\phi}\displaystyle\frac{a^{\frac{1}{3}}}{|x-x_{0}|^{\frac{2}{3}}}\gamma^{2}dxdt=\int_{Q}e^{2s\phi}\displaystyle\frac{p(x)}{(x-x_{0})^{2}}\gamma^{2}dxdt\leq C\int_{Q}e^{2s\phi}p(x)\gamma^{2}_{x}dxdt\\[5mm]&\quad\quad\quad\quad\quad\quad\quad\quad\quad\quad~~\leq CC_{1}\displaystyle\int_{Q}e^{2s\phi}a\gamma^{2}_{x}dxdt,
\end{array}
\end{eqnarray}
where $C_{1}:=\max\left\{\left(\displaystyle\frac{x^{2}_{0}}{a(0)}\right)^{\frac{2}{3}},\left(\displaystyle\frac{(1-x_{0})^{2}}{a(1)}\right)^{\frac{2}{3}}\right\}$.~Combining ~(\ref{8.24})~with~(\ref{8.25}),~we obtain
\begin{eqnarray}\label{8.26}
\displaystyle\int_{Q}e^{2s\phi}|\gamma\chi_{O_{d}}|^{2}dxdt\leq C\int_{Q}e^{2s\phi}a\gamma^{2}_{x}dxdt+C\int_{Q}e^{2s\phi}\frac{|x-x_{0}|^{2}}{a}\gamma^{2}dxdt.
\end{eqnarray}
Therefore, for a sufficiently large $s$, the right term $\displaystyle\int_{Q}e^{2s\phi}|\gamma\chi_{O_{d}}|^{2}dxdt$~in~(\ref{8.23})~can be absorbed by~$I(\gamma)$. Similarly, the last three terms on the right side of~(\ref{8.23})~can also be absorbed by ~$I(\varphi)$~and~$I(\gamma)$. We conclude that 
\begin{eqnarray}\label{8.27}
I(\varphi)+I(\gamma)\leq C\displaystyle\int_{\omega^{\prime}\times(0,T)}s^{3}\theta^{3}e^{2s\phi}\varphi^{2}dxdt+C\displaystyle\int_{\omega^{\prime}\times(0,T)}s^{3}\theta^{3}e^{2s\phi}\gamma^{2}dxdt.
\end{eqnarray}
Next, we estimate the last term in (\ref{8.27}). For this purpose, we choose a cut-off function~$\xi\in C_{0}^{\infty}(\mathbb{R})$ such that $0\leq\xi(x)\leq1,supp\xi\subseteq\omega^{\prime\prime},~\xi(x)\equiv1~in~\omega^{\prime}, ~\displaystyle\left|\frac{\xi_{x}}{\xi^{\frac{1}{2}}}\right|\leq C$ and $\displaystyle\left|\frac{\xi_{xx}}{\xi^{\frac{1}{2}}}\right|\leq C,$ where $\omega^{\prime\prime}$ is a 
 nonempty open subset of $(0,1)$ satisfying $\overline{\omega^{\prime}}\subseteq\omega^{\prime\prime},~\overline{\omega^{\prime\prime}}\subseteq O_{d}\cap\omega$ and $x_{0}\notin\overline{\omega^{\prime\prime}}$. 
 By ~(\ref{8.21}), noting that $(e^{2s\phi}s^{3}\theta^{3})|_{t=0}=(e^{2s\phi}s^{3}\theta^{3})|_{t=T}=0$~in~$(0,1)$, we have
\begin{eqnarray*}
&&0=\displaystyle\int_{Q}(\xi e^{2s\phi}s^{3}\theta^{3}\gamma\varphi)_{t}dxdt\\[2mm]
&&~~~=\displaystyle\int_{Q}\xi e^{2s\phi}s^{3}\theta^{3}(\gamma\varphi_{t}+\varphi\gamma_{t})dxdt+\int_{Q}(\xi e^{2s\phi}s^{3}\theta^{3})_{t}\gamma\varphi dxdt\\[2mm]
&&~~~=\displaystyle\int_{Q}\xi e^{2s\phi}s^{3}\theta^{3}\!\!\left\{\gamma\left[\!-\!(a\varphi_{x})_{x}\!+\!c\varphi\!-\!\gamma\chi_{O_{d}}\right]\!+\!\varphi\left[(a\gamma_{x})_{x}\!-\!c\gamma\!-\!\left(\frac{\alpha_{1}}{\mu_{1}}\chi_{\omega_{1}}\!+\!\frac{\alpha_{2}}{\mu_{2}}\chi_{\omega_{2}}\right)\varphi\right]\right\}dxdt\\[2mm]
&&~\quad+\int_{Q}(\xi e^{2s\phi}s^{3}\theta^{3})_{t}\gamma\varphi dxdt\\[2mm]
&&~~~=\int_{Q}-(\xi e^{2s\phi}s^{3}\theta^{3}\gamma a\varphi_{x})_{x}dxdt\!+\!\int_{Q}(\xi e^{2s\phi}s^{3}\theta^{3})_{x}\gamma a\varphi_{x}dxdt\!+\!\int_{Q}\xi e^{2s\phi}s^{3}\theta^{3}\gamma_{x}a\varphi_{x}dxdt\\[2mm]
&&~\quad+\int_{Q}(\xi e^{2s\phi}s^{3}\theta^{3}\varphi a\gamma_{x})_{x}dxdt\!-\!\int_{Q}(\xi e^{2s\phi}s^{3}\theta^{3})_{x}\varphi a\gamma_{x}dxdt\!-\!\int_{Q}\xi e^{2s\phi}s^{3}\theta^{3}\varphi_{x}a\gamma_{x}dxdt\\[2mm]
&&~\quad-\int_{Q}\xi e^{2s\phi}s^{3}\theta^{3}\gamma^{2}\chi_{O_{d}}dxdt\!-\!\int_{Q}\xi e^{2s\phi}s^{3}\theta^{3}\left(\frac{\alpha_{1}}{\mu_{1}}\chi_{\omega_{1}}+\frac{\alpha_{2}}{\mu_{2}}\chi_{\omega_{2}}\right)\varphi^{2}dxdt\\[2mm]
&&~\quad+\int_{Q}(\xi e^{2s\phi}s^{3}\theta^{3})_{t}\gamma\varphi dxdt\\[2mm]
&&~~~=
\int_{Q}\!-\!(\xi e^{2s\phi}s^{3}\theta^{3})_{xx}\gamma a\varphi dxdt\!\!-\!\!2\int_{Q}(\xi e^{2s\phi}s^{3}\theta^{3})_{x}\gamma_{x}a\varphi dxdt\!\!-\!\!\int_{Q}(\xi e^{2s\phi}s^{3}\theta^{3})_{x}\gamma a_{x}\varphi dxdt\\[2mm]
&&~\quad-\int_{Q}\xi e^{2s\phi}s^{3}\theta^{3}\gamma^{2}\chi_{O_{d}}dxdt+\int_{Q}(\xi e^{2s\phi}s^{3}\theta^{3})_{t}\gamma\varphi dxdt\\[2mm]
&&~\quad-\int_{Q}\xi e^{2s\phi}s^{3}\theta^{3}\left(\frac{\alpha_{1}}{\mu_{1}}\chi_{\omega_{1}}+\frac{\alpha_{2}}{\mu_{2}}\chi_{\omega_{2}}\right)\varphi^{2}dxdt.
\end{eqnarray*}
This implies
\begin{eqnarray*}
&&\displaystyle\int_{Q}\xi e^{2s\phi}s^{3}\theta^{3}\left[\gamma^{2}\chi_{O_{d}}+\left(\frac{\alpha_{1}}{\mu_{1}}\chi_{\omega_{1}}+\frac{\alpha_{2}}{\mu_{2}}\chi_{\omega_{2}}\right)\varphi^{2}\right]dxdt\\[2mm]
&&=\displaystyle\int_{Q}-(\xi e^{2s\phi}s^{3}\theta^{3})_{xx}\gamma a\varphi dxdt-2\int_{Q}(\xi e^{2s\phi}s^{3}\theta^{3})_{x}\gamma_{x}a\varphi dxdt\\[2mm]
&&\quad-\displaystyle\int_{Q}(\xi e^{2s\phi}s^{3}\theta^{3})_{x}\gamma a_{x}\varphi dxdt+\int_{Q}(\xi e^{2s\phi}s^{3}\theta^{3})_{t}\gamma\varphi dxdt\\[2mm]
&&=-\displaystyle\int_{Q}\gamma a\varphi\left[\xi(e^{2s\phi}s^{3}\theta^{3})_{xx}+2\xi_{x}(e^{2s\phi}s^{3}\theta^{3})_{x}+\xi_{xx}e^{2s\phi}s^{3}\theta^{3}\right]dxdt\\[2mm]
&&~~~~-2\displaystyle\int_{Q}\gamma_{x}a\varphi\left[\xi(e^{2s\phi}s^{3}\theta^{3})_{x}+\xi_{x}e^{2s\phi}s^{3}\theta^{3}\right]dxdt
-\displaystyle\int_{Q}\gamma a_{x}\varphi\left[\xi(e^{2s\phi}s^{3}\theta^{3})_{x}\right.\\[2mm]
&&~~~~+\displaystyle\left.\xi_{x}e^{2s\phi}s^{3}\theta^{3}\right]dxdt+\displaystyle\int_{Q}\xi(e^{2s\phi}s^{3}\theta^{3})_{t}\gamma \varphi dxdt.
\end{eqnarray*}
Further,
\begin{eqnarray}\label{8.28}
\begin{array}{rl}
&\displaystyle\int_{Q}\xi e^{2s\phi}s^{3}\theta^{3}\gamma^{2}\chi_{O_{d}}dxdt\\[5mm]
&\leq\displaystyle-\int_{Q}\xi\left[(e^{2s\phi}s^{3}\theta^{3})_{xx}\gamma a\varphi+2(e^{2s\phi}s^{3}\theta^{3})_{x}\gamma_{x}a\varphi+(e^{2s\phi}s^{3}\theta^{3})_{x}\gamma a_{x}\varphi\right]dxdt\\[5mm]
&~~~~-\displaystyle\int_{Q}\xi_{x}\left[2(e^{2s\phi}s^{3}\theta^{3})_{x}\gamma a\varphi+e^{2s\phi}s^{3}\theta^{3}(2\gamma_{x}a\varphi+\gamma a_{x}\varphi)\right]dxdt\\[5mm]
&~~~~-\displaystyle\int_{Q}\xi_{xx}e^{2s\phi}s^{3}\theta^{3}\gamma a\varphi dxdt+\displaystyle\int_{Q}\xi(e^{2s\phi}s^{3}\theta^{3})_{t}\gamma \varphi dxdt.
\end{array}
\end{eqnarray}
Notice that $x_{0}\notin\overline{\omega^{\prime\prime}}$, it is easy to check that there exists a positive constant~$C$~such that the following estimate holds in ~$\omega^{\prime\prime}\times(0,T)$
\begin{equation}\label{8.29}
|(e^{2s\phi}s^{3}\theta^{3})_{t}|\leq Ce^{2s\phi}s^{4}\theta^{5},~~
|(e^{2s\phi}s^{3}\theta^{3})_{x}|\leq Ce^{2s\phi}s^{4}\theta^{4},~~
|(e^{2s\phi}s^{3}\theta^{3})_{xx}|\leq Ce^{2s\phi}s^{5}\theta^{5}.
\end{equation}
By~(\ref{8.28})~and~(\ref{8.29})~we obtain
\begin{eqnarray}\label{8.30}
\begin{array}{rl}
&\displaystyle\int_{Q}\xi e^{2s\phi}s^{3}\theta^{3}\gamma^{2}\chi_{O_{d}}dxdt
\\[5mm]&\leq C\displaystyle\int_{Q}e^{2s\phi}\xi\left(s^{5}\theta^{5}|\gamma a\varphi|+s^{4}\theta^{4}|\gamma_{x}a\varphi|+s^{4}\theta^{4}|\gamma a_{x}\varphi|\right)dxdt\\[5mm]
&~~~~+C\displaystyle\int_{Q}e^{2s\phi}\Big{|}\frac{\xi_{x}}{\xi^{\frac{1}{2}}}\Big{|}\xi^{\frac{1}{2}}\left(s^{4}\theta^{4}|\gamma a\varphi|+s^{3}\theta^{3}|\gamma_{x}a\varphi|+s^{3}\theta^{3}|\gamma a_{x}\varphi|\right)dxdt\\[5mm]
&~~~~+C\displaystyle\int_{Q}e^{2s\phi}\left(\Big{|}\frac{\xi_{xx}}{\xi^{\frac{1}{2}}}\Big{|}\xi^{\frac{1}{2}}s^{3}\theta^{3}|\gamma a\varphi| +\xi s^{4}\theta^{5}|\gamma\varphi|\right)dxdt\\[5mm]
&\leq C\displaystyle\int_{\omega^{\prime\prime}\times(0,T)}e^{2s\phi}\left(s^{5}\theta^{5}|\gamma\varphi|+s^{4}\theta^{4}|\gamma_{x}\sqrt{a}\varphi|+s^{4}\theta^{4}|\gamma\varphi|\right)dxdt
\\[5mm]&\quad+C\displaystyle\int_{\omega^{\prime\prime}\times(0,T)}e^{2s\phi}\left(s^{4}\theta^{4}|\gamma \varphi|+s^{3}\theta^{3}|\gamma_{x}\sqrt{a}\varphi|+s^{3}\theta^{3}|\gamma\varphi|\right)dxdt
\\[5mm]&\quad+C\displaystyle\int_{\omega^{\prime\prime}\times(0,T)}e^{2s\phi}s^{3}\theta^{3}|\gamma \varphi|dxdt+C\displaystyle\int_{\omega^{\prime\prime}\times(0,T)}e^{2s\phi}s^{4}\theta^{5}|\gamma\varphi|dxdt
\\[5mm]
&\leq C\displaystyle\int_{\omega^{\prime\prime}\times(0,T)}e^{2s\phi}\left(s^{5}\theta^{5}|\gamma\varphi|+s^{4}\theta^{4}|\gamma_{x}\sqrt{a}\varphi|
\right)dxdt\\[5mm]&:=I_{1}+I_{2}.
\end{array}
\end{eqnarray}
Applying Cauchy inequality with~$\varepsilon$, we get
\begin{eqnarray}\label{8.31}
I_{1}\leq \displaystyle\int_{\omega^{\prime\prime}\times(0,T)}e^{2s\phi}\left(\frac{1}{2\varepsilon}s^{7}\theta^{7}\varphi^{2}+\frac{\varepsilon}{2}s^{3}\theta^{3}\gamma^{2}\right)dxdt,
\end{eqnarray}
and
\begin{eqnarray}\label{8.32}
I_{2}\leq \displaystyle\int_{\omega^{\prime\prime}\times(0,T)}e^{2s\phi}\left(\frac{1}{2\varepsilon}s^{7}\theta^{7}\varphi^{2}+\frac{\varepsilon}{2}s\theta a\gamma_{x}^{2}\right)dxdt.
\end{eqnarray}
Substituting~(\ref{8.31})~and~(\ref{8.32})~into~(\ref{8.30})~and taking~$\varepsilon$~sufficiently small, we obtain
\begin{eqnarray}\label{8.33}
\begin{array}{rl}
&\displaystyle\int_{Q}\xi e^{2s\phi}s^{3}\theta^{3}\gamma^{2}\chi_{O_{d}}dxdt\\[5mm]&\leq C\left(\displaystyle\int_{\omega^{\prime\prime}\times(0,T)}s^{7}\theta^{7}e^{2s\phi}\varphi^{2}dxdt+\varepsilon\displaystyle\int_{\omega^{\prime\prime}\times(0,T)}s\theta e^{2s\phi}a\gamma_{x}^{2}dxdt\right).
\end{array}
\end{eqnarray}
Substituting (\ref{8.33})~into~(\ref{8.27}), we have
\begin{eqnarray}\label{8.34}
\begin{array}{rl}
I(\varphi)+I(\gamma)
 \leq C\displaystyle\int_{\omega\times(0,T)}s^{7}\theta^{7}e^{2s\phi}\varphi^{2}dxdt
.
\end{array}
\end{eqnarray}
This completes the proof.
\endpf

The main result of this section is the following observability inequality for the system (\ref{8.16}).

\begin{theorem}\label{ttt}
Assume that $O_{d}\cap\omega\neq\emptyset,~x_{0}\in\omega$ and Condition~3~holds. If
 $\mu_{i}(i=1,2)$ are sufficiently large, then estimate (\ref{8.17}) holds for (\ref{8.16}).
\end{theorem}

\no{\bf Proof.} We introduce a cut-off
function~$\delta(\cdot)\in C^{1}[0,T]$ such that
\begin{eqnarray}\label{8.35}
0\leq\delta(t)\leq1~\mbox{in}~[0,T];~\delta(t)=1~\mbox{in}~[0,\frac{T}{2}];~\delta(t)=0~\mbox{in}~[\frac{3T}{4},T];~|\delta_{t}(t)|\leq\frac{C}{T}.
\end{eqnarray}
Multiplying both sides of the first equation of (\ref{8.21}) by~$\delta\varphi$~and integrating it on $(0,1)$, we find that
\begin{eqnarray*}
\int_{0}^{1}(-\varphi_{t}-(a\varphi_{x})_{x}+c\varphi)\delta\varphi dx=\int_{0}^{1}\gamma\chi_{O_{d}}\delta\varphi dx.
\end{eqnarray*}
By Cauchy inequality, it holds that 
\begin{eqnarray}\label{8.37}
\begin{array}{rl}
-\displaystyle\frac{d}{dt}\displaystyle\int_{0}^{1}\delta\varphi^{2}dx+2\displaystyle\int_{0}^{1}\delta a\varphi^{2}_{x}dx\leq-\delta_{t}\displaystyle\int_{0}^{1}\varphi^{2}dx+(1+2\|c\|_{{L^{\infty}(Q)}})\displaystyle\int_{0}^{1}\delta\varphi^{2}dx+\displaystyle\int_{0}^{1}\delta\gamma^{2}dx.
\end{array}
\end{eqnarray}
Multiplying both sides of~(\ref{8.37})~by~$e^{(1+2\|c\|_{L^{\infty}(Q)})t}$~and integrating it over~$[0,T]$,~we obtain
\begin{eqnarray*}
\begin{array}{rl}
&\!\!\displaystyle\int_{0}^{1}|\varphi(x,0)|^{2}dx\!\!+\!\!2\!\!\displaystyle\int_{(0,1)\times(0,\frac{T}{2})}e^{(1+2\|c\|_{L^{\infty}(Q)})t}a\varphi^{2}_{x}dxdt\!\!+\!\!2\!\!\displaystyle\int_{(0,1)\times(\frac{T}{2},\frac{3T}{4})}e^{(1+2\|c\|_{L^{\infty}(Q)})t}\delta a\varphi^{2}_{x}dxdt\\[5mm]&\quad\leq C\left(\displaystyle\int_{(0,1)\times(\frac{T}{2},\frac{3T}{4})}e^{(1+2\|c\|_{L^{\infty}(Q)})t}\varphi^{2}dxdt+\displaystyle\int_{(0,1)\times(0,\frac{3T}{4})}e^{(1+2\|c\|_{L^{\infty}(Q)})t}\gamma^{2}dxdt\right).
\end{array}
\end{eqnarray*}
Notice that $e^{(1+2\|c\|_{L^{\infty}(Q)})t}$~is bounded in~$[0,T]$, we get
\begin{eqnarray}\label{8.38}
\!\!\int_{0}^{1}|\varphi(x,0)|^{2}dx\!\!+\!\!2\int_{(0,1)\times(0,\frac{T}{2})}\!a\varphi^{2}_{x}dxdt\!\!\leq\!\! C\left(\int_{(0,1)\times(\frac{T}{2},\frac{3T}{4})}\varphi^{2}dxdt\!\!+\!\!\int_{(0,1)\times(0,\frac{3T}{4})}\gamma^{2}dxdt\right).
\end{eqnarray}
Since $\left(\displaystyle\frac{a}{|x-x_{0}|^{2}}\right)^{\frac{1}{3}}$~has a positive lower bound, we take~$p(x)=a(x)\left(\displaystyle\frac{|x-x_{0}|^{2}}{a}\right)^{\frac{2}{3}}$, similar to (\ref{8.25}), by Lemma~\ref{l2} we conclude that 
\begin{eqnarray}\label{8.39}
\begin{array}{rl}
&\displaystyle\int_{(0,1)\times(0,\frac{T}{2})}\displaystyle\frac{(x-x_{0})^{2}}{a}\varphi^{2}dxdt\leq C\displaystyle\int_{(0,1)\times(0,\frac{T}{2})}\varphi^{2}dxdt\\[5mm]&\leq C\displaystyle\int_{(0,1)\times(0,\frac{T}{2})}\left(\displaystyle\frac{a}{|x-x_{0}|^{2}}\right)^{\frac{1}{3}}\varphi^{2}dxdt=C\displaystyle\int_{(0,1)\times(0,\frac{T}{2})}\displaystyle\frac{p(x)} {|x-x_{0}|^{2}}\varphi^{2}dxdt\\[5mm]&\leq C\displaystyle\int_{(0,1)\times(0,\frac{T}{2})}p(x)\varphi_{x}^{2}dxdt\leq C\displaystyle\int_{(0,1)\times(0,\frac{T}{2})}a\varphi_{x}^{2}dxdt.
\end{array}
\end{eqnarray}
Combining~(\ref{8.38})~with~(\ref{8.39}), we get
\begin{eqnarray}\label{8.40}
\begin{array}{rl}
&\displaystyle\int_{0}^{1}|\varphi(x,0)|^{2}dx+\displaystyle\int_{(0,1)\times(0,\frac{T}{2})}a\varphi^{2}_{x}dxdt+\displaystyle\int_{(0,1)\times(0,\frac{T}{2})}\displaystyle\frac{(x-x_{0})^{2}}{a}\varphi^{2}dxdt\\[5mm]&\leq C\left(\displaystyle\int_{(0,1)\times(\frac{T}{2},\frac{3T}{4})}\varphi^{2}dxdt+\displaystyle\int_{(0,1)\times(0,\frac{3T}{4})}\gamma^{2}dxdt\right).
\end{array}
\end{eqnarray}

Let~$\eta(\cdot)\in C^{1}[0,T]$~be a function satisfying
\begin{eqnarray}\label{8.36}
\eta(t)=\left\{\begin{array}{ll}
\displaystyle\frac{2^{8}}{T^{8}},&0\leq t\leq\frac{T}{2},\\[2mm]
\displaystyle\frac{1}{t^{4}(T-t)^{4}},&\frac{T}{2}\leq t\leq T.
\end{array}\right.
\end{eqnarray}
We define a new weight function as follows,
$$\bar{\phi}(x,t)=\eta(t)\psi(x),$$
where~$\psi(x)$ is given in (\ref{7.2}). Since~$e^{2s\bar{\phi}}$~has a positive upper bound in~$(0,1)\times(0,\frac{T}{2})$, it follows from~(\ref{8.40})~that
\begin{eqnarray}\label{8.41}
\begin{array}{rl}
&\displaystyle\int_{0}^{1}|\varphi(x,0)|^{2}dx+\displaystyle\int_{(0,1)\times(0,\frac{T}{2})}e^{2s\bar{\phi}}\left(s^{3}\eta^{3}\displaystyle\frac{(x-x_{0})^{2}}{a}\varphi^{2}+s\eta a\varphi^{2}_{x}\right)dxdt\\[5mm]&\leq C\left(\displaystyle\int_{(0,1)\times(\frac{T}{2},\frac{3T}{4})}\varphi^{2}dxdt+\displaystyle\int_{(0,1)\times(0,\frac{3T}{4})}\gamma^{2}dxdt\right).
\end{array}
\end{eqnarray}
Define~ $\bar{I}_{[m,n]}(z)=\displaystyle\int_{(0,1)\times(m,n)}e^{2s\bar{\phi}}\left(s^{3}\eta^{3}\displaystyle\frac{(x-x_{0})^{2}}{a}z^{2}+s\eta az^{2}_{x}\right)dxdt.$
Similar to~(\ref{8.24})~and~(\ref{8.25}), note that $e^{2s\bar{\phi}}$~has a positive lower bound in~$(0,1)\times(0,\frac{T}{2})$, by Young~inequality and Lemma~\ref{l2},  we arrive at 
\begin{eqnarray}\label{8.42}
\begin{array}{rl}
&\displaystyle\int_{(0,1)\times(0,\frac{T}{2})}\gamma^{2}dxdt\leq\displaystyle\frac{3}{4}{C}_{1}\displaystyle\int_{(0,1)\times(0,\frac{T}{2})}a\gamma^{2}_{x}dxdt+\displaystyle\frac{1}{4}\displaystyle\int_{(0,1)\times(0,\frac{T}{2})}\frac{(x-x_{0})^{2}}{a}\gamma^{2}dxdt
\\[5mm]&\quad\quad\quad\quad\quad\quad\quad\quad\leq C\bar{I}_{[0,\frac{T}{2}]}(\gamma),
\end{array}
\end{eqnarray}
where~$C_{1}:=\max\left\{\left(\displaystyle\frac{x^{2}_{0}}{a(0)}\right)^{\frac{2}{3}},\left(\displaystyle\frac{(1-x_{0})^{2}}{a(1)}\right)^{\frac{2}{3}}\right\}$.
 We add~$\bar{I}_{[0,\frac{T}{2}]}(\gamma)$~to both sides of the inequality (\ref{8.41}), by  (\ref{8.42}) and $e^{2s\bar{\phi}}$~has positive lower bound in~$(0,1)\times(0,\frac{T}{2})$, we get
 \begin{eqnarray}\label{8.43}
\begin{array}{rl}
&\displaystyle\int_0^{1}|\varphi(x,0)|^{2}dx+\bar{I}_{[0,\frac{T}{2}]}(\varphi)+\bar{I}_{[0,\frac{T}{2}]}(\gamma)\\[5mm]&\leq \bar{I}_{[0,\frac{T}{2}]}(\gamma)+C\left[\displaystyle\int_{(0,1)\times(\frac{T}{2},\frac{3T}{4})}(\varphi^{2}+\gamma^{2})dxdt+\displaystyle\int_{(0,1)\times(0,\frac{T}{2})}\gamma^{2}dxdt\right]
\\[5mm]&\leq C\bar{I}_{[0,\frac{T}{2}]}(\gamma)
+C\displaystyle\int_{(0,1)\times(\frac{T}{2},\frac{3T}{4})}(\varphi^{2}+\gamma^{2})dxdt.
\end{array}
\end{eqnarray}

On the other hand, since~$e^{2s\bar{\phi}}$~has a positive upper bound in~$(0,1)\times(0,\frac{T}{2})$, by the energy estimate for the solution of second equation in system (\ref{8.21}), it is easy to obtain that 
\begin{eqnarray}\label{8.44}
\bar{I}_{[0,\frac{T}{2}]}(\gamma)\leq C\int_{(0,1)\times(0,\frac{T}{2})}(\gamma^{2}+a\gamma^{2}_{x})dxdt\leq C\left(\frac{\alpha^{2}_{1}}{\mu^{2}_{1}}+\frac{\alpha^{2}_{2}}{\mu^{2}_{2}}\right)\int_{(0,1)\times(0,\frac{T}{2})}\sigma^{-4}\varphi^{2}dxdt.
\end{eqnarray}
Similar to~(\ref{8.42}), notice that $\sigma_{0}=\displaystyle\min_{t\in[0,T]}\sigma(t)$~we get
\begin{eqnarray}\label{8.45}
\begin{array}{rl}
&\left(\displaystyle\frac{\alpha^{2}_{1}}{\mu^{2}_{1}}+\frac{\alpha^{2}_{2}}{\mu^{2}_{2}}\right)\displaystyle\int_{(0,1)\times(0,\frac{T}{2})}\sigma^{-4}\varphi^{2}dxdt\\[5mm]&\leq\left(\displaystyle\frac{\alpha^{2}_{1}}{\mu^{2}_{1}}+\frac{\alpha^{2}_{2}}{\mu^{2}_{2}}\right)\sigma_{0}^{-4}\left(\displaystyle\frac{3}{4}C_{1}\displaystyle\int_{(0,1)\times(0,\frac{T}{2})}a\varphi^{2}_{x}dxdt+\displaystyle\frac{1}{4}\displaystyle\int_{(0,1)\times(0,\frac{T}{2})}\frac{(x-x_{0})^{2}}{a}\varphi^{2}dxdt\right)\\[5mm]&\leq C\left(\displaystyle\frac{\alpha^{2}_{1}}{\mu^{2}_{1}}+\frac{\alpha^{2}_{2}}{\mu^{2}_{2}}\right)\bar{I}_{[0,\frac{T}{2}]}(\varphi).
\end{array}
\end{eqnarray}
By (\ref{8.43})-(\ref{8.45}), we choose~$\mu_{i}(i=1,2)$ large enough, then \begin{eqnarray}\label{8.46}
\begin{array}{rl}
&\displaystyle\int_0^{1}|\varphi(x,0)|^{2}dx+\bar{I}_{[0,\frac{T}{2}]}(\varphi)+\bar{I}_{[0,\frac{T}{2}]}(\gamma)\leq
C\displaystyle\int_{(0,1)\times(\frac{T}{2},\frac{3T}{4})}(\varphi^{2}+\gamma^{2})dxdt.
\end{array}
\end{eqnarray}
Since~$e^{2s\phi}$~has a positive lower bound, similar to~(\ref{8.45}), we obtain
\begin{eqnarray}\label{8.47}
\begin{array}{rl}
&\displaystyle\int_{(0,1)\times(\frac{T}{2},\frac{3T}{4})}(\varphi^{2}+\gamma^{2})dxdt\\[5mm]&\leq\displaystyle\frac{3}{4}C_{1}\displaystyle\int_{(0,1)\times(\frac{T}{2},\frac{3T}{4})}a\varphi^{2}_{x}dxdt+\displaystyle\frac{1}{4}\displaystyle\int_{(0,1)\times(\frac{T}{2},\frac{3T}{4})}\frac{(x-x_{0})^{2}}{a}\varphi^{2}dxdt\\[5mm]&~~+\displaystyle\frac{3}{4}C_{1}\displaystyle\int_{(0,1)\times(\frac{T}{2},\frac{3T}{4})}a\gamma^{2}_{x}dxdt+\displaystyle\frac{1}{4}\displaystyle\int_{(0,1)\times(\frac{T}{2},\frac{3T}{4})}\frac{(x-x_{0})^{2}}{a}\gamma^{2}dxdt\\[5mm]&\leq C\displaystyle\int_{(0,1)\times(\frac{T}{2},\frac{3T}{4})}\left[s\theta a\varphi^{2}_{x}+s^{3}\theta^{3}\frac{(x-x_{0})^{2}}{a}\varphi^{2}+s\theta a\gamma^{2}_{x}+s^{3}\theta^{3}\frac{(x-x_{0})^{2}}{a}\gamma^{2}\right]e^{2s\phi}dxdt\\[5mm]&\leq C(I(\varphi)+I(\gamma)).\end{array}
\end{eqnarray}
By (\ref{8.46})~and~(\ref{8.47}), together with Proposition~\ref{p3}, indicates
\begin{eqnarray}\label{8.48}
\begin{array}{rl}
&\displaystyle\int_{0}^{1}|\varphi(x,0)|^{2}dx+\bar{I}_{[0,\frac{T}{2}]}(\varphi)+\bar{I}_{[0,\frac{T}{2}]}(\gamma) \leq C\displaystyle\int_{\omega\times(0,T)}\varphi^{2}dxdt.
\end{array}
\end{eqnarray}

Since~$\eta=\theta,\bar{\phi}=\phi$~in~$(0,1)\times(\displaystyle\frac{T}{2},T)$, by Proposition~\ref{p3} again, it follows that 
\begin{eqnarray}\label{8.49}
\bar{I}_{[\frac{T}{2},T]}(\varphi)+\bar{I}_{[\frac{T}{2},T]}(\gamma)\leq I(\varphi)+I(\gamma)
\leq C\int_{\omega\times(0,T)}\varphi^{2}dxdt.
\end{eqnarray}
Adding (\ref{8.48}) and (\ref{8.49}), we obtain
\begin{eqnarray}\label{8.50}
\int_{0}^{1}|\varphi(x,0)|^{2}dx+\bar{I}_{[0,T]}(\varphi)+\bar{I}_{[0,T]}(\gamma)
\leq C\int_{\omega\times(0,T)}\varphi^{2}dxdt.
\end{eqnarray}

Finally, it is easy to obtain the following weighted energy estimates for the solutions~$\psi^{i}(i=1,2)$~of~(\ref{8.16})
\begin{eqnarray}\label{8.51}
&\displaystyle\int_{Q}|\psi^{i}|^{2}dxdt\leq C\displaystyle\int_{Q}\sigma^{-4}\varphi^{2}dxdt.
\end{eqnarray}
Similar to~(\ref{8.26}), notice that $\phi^{*}\leq\phi\leq\bar{\phi}$, we get
\begin{eqnarray}\label{8.52}
\begin{array}{rl}
&\displaystyle\int_{Q}\sigma^{-4}\varphi^{2}dxdt\leq C\displaystyle\int_{Q}e^{2s\phi^{*}}\left(a\varphi^{2}_{x}+\frac{|x-x_{0}|^{2}}{a}\varphi^{2}\right)dxdt\leq C\bar{I}_{[0,T]}(\varphi).
\end{array}
\end{eqnarray}
By~(\ref{8.50})-(\ref{8.52}), we have
\begin{eqnarray}\label{8.53}
\begin{array}{rl}
&\displaystyle\int_{Q}|\psi^{i}|^{2}dxdt\leq C\bar{I}_{[0,T]}(\varphi)\leq C\int_{\omega\times(0,T)}\varphi^{2}dxdt.
\end{array}
\end{eqnarray}
Combining (\ref{8.50})~and~(\ref{8.53}), the proof is completed.
\endpf
\medskip

\section{Null controllability of the Semilinear system}

The main purpose of this section is to prove Theorem \ref{T1} for  the semilinear case, (i.e., $G(y)\neq0$). 

\no{\bf Proof of Theorem \ref{T1} for the semilinear case.} By Definition~\ref{d2}, similar to~$(\ref{8.3})$, we can see that for any given leader control~$f$,~$(\tilde{u}_{1}(f),\tilde{u}_{2}(f))$~is the Nash quasi equilibrium pair of~$J_{i}$ if and only if 
\begin{eqnarray}\label{8.54}
\mu_{i}\int_{\omega_{i}\times(0,T)}\sigma^{2}\tilde{u}_{i} \hat{u}_{i}dxdt+\alpha_{i}\int_{O_{i,d}\times(0,T)}(y-y_{i,d}) \hat{y}^{i}dxdt=0,  ~\forall\hat{u}_{i}\in S_{i},~i=1,2,
\end{eqnarray}
\noindent
where~$\hat{y}^{i}$~denotes the solution of the following system corresponding to~$\hat{u}_{i}$:
\begin{equation}\label{8.55}
\left\{\begin{array}{ll}
\hat{y}^{i}_{t}-(a\hat{y}^{i}_{x})_{x}+c\hat{y}^{i}=G^{\prime}(y)\hat{y}^{i}+\hat{u}_{i}\chi_{\omega_{i}}& \mbox{ in }Q,\\[3mm]
\hat{y}^{i}(0,t)=\hat{y}^{i}(1,t)=0&\mbox{ in }(0,T),\\[3mm]
\hat{y}^{i}(x,0)=0&\mbox{ in }(0,1),
\end{array}\right.
\end{equation}
\noindent
and $y$~is the solution of system~(\ref{6.1})~corresponding to~$(f,\tilde{u}_{1}(f),\tilde{u}_{2}(f))$. 
We introduce the following adjoint of system (\ref{8.55}):
\begin{equation}\label{8.57}
\left\{\begin{array}{ll}
-\hat{q}^{i}_{t}-(a\hat{q}^{i}_{x})_{x}+c\hat{q}^{i}=G^{\prime}(y)\hat{q}^{i}+\alpha_{i}(y-y_{i,d})\chi_{O_{i,d}}& \mbox{ in }Q,\\[3mm]
\hat{q}^{i}(0,t)=\hat{q}^{i}(1,t)=0&\mbox{ in }(0,T),\\[3mm]
\hat{q}^{i}(x,T)=0&\mbox{ in }(0,1),
\end{array}\right.
\end{equation}
\noindent
where~$y$~is the solution of the system~(\ref{6.1})~corresponding to~$(f,\tilde{u}_{1}(f),\tilde{u}_{2}(f))$.
Multiplying both sides of the first equation of~(\ref{8.55})~by~$\hat{q}^{i}$~and integrating it over~$Q$,~we obtain
\begin{eqnarray}\label{8.56}
\begin{array}{rl}
&\displaystyle\int_Q(-\hat{q}^{i}_{t}-(a\hat{q}^{i}_{x})_{x}+c\hat{q}^{i})\hat{y}^{i}dxdt+\int_{0}^{1}\hat{y}^{i}\hat{q}^{i}dx \bigg|_{0}^{T}-\displaystyle\int_{0}^{T}\hat{q}^{i}\cdot a\hat{y}^{i}_{x}dt\bigg|_{0}^{1}+\displaystyle\int_{0}^{T}\hat{y}^{i}\cdot a\hat{q}^{i}_{x}dt\bigg|_{0}^{1}\\[5mm]
&=\displaystyle\int_QG^{\prime}(y)\hat{y}^{i}\hat{q}^{i}dxdt+\int_{\omega_{i}\times(0,T)}\hat{u}_{i}\hat{q}^{i}dxdt.
\end{array}
\end{eqnarray}
By (\ref{8.54})~and equation~(\ref{8.55}), (\ref{8.56})~can be reduced to
$$\int_{\omega_{i}\times(0,T)}(\hat{q}^{i}+\mu_{i}\sigma^{2}\tilde{u}_{i})\hat{u}_{i}dxdt=0,~\forall \hat{u}_{i}\in S_{i},~i=1,2.$$
Therefore, the explicit expression of~$\tilde{u}_{i}$~is as follows
\begin{eqnarray}\label{8.58}
\tilde{u}_{i}=-\frac{1}{\mu_{i}}\sigma^{-2}\hat{q}^{i}\chi_{\omega_{i}},~i=1,2.
\end{eqnarray}
As a conclusion, we obtain the following optimality system
\begin{equation}\label{8.59}
\left\{\begin{array}{ll}
y_{t}-(ay_{x})_{x}+cy=G(y)+f\chi_{\omega}-\frac{1}{\mu_{1}}\sigma^{-2}\hat{q}^{1}\chi_{\omega_{1}}-\frac{1}{\mu_{2}}\sigma^{-2}\hat{q}^{2}\chi_{\omega_{2}}& \mbox{ in }Q,\\[3mm]
-\hat{q}^{i}_{t}-(a\hat{q}^{i}_{x})_{x}+c\hat{q}^{i}=G^{\prime}(y)\hat{q}^{i}+\alpha_{i}(y-y_{i,d})\chi_{O_{i,d}}& \mbox{ in }Q,\\[3mm]
y(0,t)=y(1,t)=0,~\hat{q}^{i}(0,t)=\hat{q}^{i}(1,t)=0& \mbox{ in }(0,T),\\[3mm]
y(x,0)=y_{0}(x),~\hat{q}^{i}(x,T)=0& \mbox{ in }(0,1).
\end{array}\right.
\end{equation}

In order to prove the null controllability of~(\ref{6.1}), we only need to prove that the optimal system~(\ref{8.59})~is null controllable with respect to the solution $y$, i.e., our objective is to find a control~$\tilde{f}\in L^{2}(\omega\times(0,T))$~such that the solution of~(\ref{8.59})  satisfies
\begin{eqnarray}\label{8.60}
y(\cdot,T;\tilde{f})=0 \mbox{~in~}(0,1),
\end{eqnarray}
and 
\begin{eqnarray}\label{8.888}
J(\tilde{f})=\displaystyle\min_{f}J(f),~~\forall f\in L^{2}(\omega\times(0,T)).
\end{eqnarray}
A fixed-point method will be used to prove this.
For any~$y\in L^{2}(Q)$~and~$f\in L^{2}(\omega\times(0,T))$, we consider the following linear system
\begin{equation}\label{8.61}
\left\{\begin{array}{ll}
\tilde{z}_{t}-(a\tilde{z}_{x})_{x}+c\tilde{z}=\tilde{G}(x,t;y)\tilde{z}+f\chi_{\omega}-\frac{1}{\mu_{1}}\sigma^{-2}\hat{q}^{1}\chi_{\omega_{1}}-\frac{1}{\mu_{2}}\sigma^{-2}\hat{q}^{2}\chi_{\omega_{2}}& \mbox{ in }Q,\\[3mm]
-\hat{q}^{i}_{t}-(a\hat{q}^{i}_{x})_{x}+c\hat{q}^{i}=G^{\prime}(y)\hat{q}^{i}+\alpha_{i}(\tilde{z}-y_{i,d})\chi_{O_{i,d}}& \mbox{ in }Q,\\[3mm]
\tilde{z}(0,t)=\tilde{z}(1,t)=0,~\hat{q}^{i}(0,t)=\hat{q}^{i}(1,t)=0& \mbox{ in }(0,T),\\[3mm]
\tilde{z}(x,0)=y_{0}(x),~\hat{q}^{i}(x,T)=0& \mbox{ in }(0,1),
\end{array}\right.
\end{equation}
where~$\tilde{G}(x,t;y)=\int_{0}^{1}G^{\prime}(\tau y)d\tau.$~Next we will study the null controllability of the solution~$\tilde{z}$~ in~(\ref{8.61}). Since~$G\in W^{1,\infty}(\mathbb{R})$, there exists a positive constant~$M$~such that
\begin{eqnarray}\label{8.62}
|\tilde{G}(x,t;\alpha)|+|G^{\prime}(\alpha)|\leq M,~\forall(x,t,\alpha)\in Q\times\mathbb{R}.
\end{eqnarray}
By the $L^2$-estimate  for parabolic equations (\ref{8.61}), we obtian
\begin{eqnarray}\label{8.63}
\begin{array}{rl}
&\displaystyle\|\tilde{z}\|_{L^{2}(0,T;H_{a}^{1}(0,1))}+\|\tilde{z}_{t}\|_{L^{2}(0,T;H^{*}_{a}(0,1))}\\[2mm]&\leq C\left(\|y_{0}\|_{L^{2}(0,1)}+\|f\|_{L^{2}(\omega\times(0,T))}+\Big{\|}\frac{1}{\mu_{1}}\hat{q}^{1}\Big{\|}_{S_{1}}+\Big{\|}\frac{1}{\mu_{2}}\hat{q}^{2}\Big{\|}_{S_{2}}\right).
\end{array}
\end{eqnarray}
By (\ref{8.58})~and~(\ref{8.8}), we conclude that, for $\mu_{1}$ and $\mu_{2}$ sufficiently large, it holds that 
\begin{eqnarray}\label{8.64}
\displaystyle\|\tilde{z}\|_{L^{2}(0,T;H_{a}^{1}(0,1))}+\|\tilde{z}_{t}\|_{L^{2}(0,T;H^{*}_{a}(0,1))}\leq C\left(\|y_{0}\|_{L^{2}(0,1)}+\|f\|_{L^{2}(\omega\times(0,T))}\right).
\end{eqnarray}

Let us introduce  the following adjoint system of~(\ref{8.61}):
\begin{equation}\label{8.65}
\left\{\begin{array}{ll}
-\tilde{\psi}_{t}-(a\tilde{\psi}_{x})_{x}+c\tilde{\psi}=\tilde{G}(x,t;y)\tilde{\psi}+\displaystyle\sum_{i=1}^{2}\alpha_{i}\tilde{\gamma}^{i}\chi_{O_{i,d}}& \mbox{ in }Q,\\[3mm]
\tilde{\gamma}^{i}_{t}-(a\tilde{\gamma}^{i}_{x})_{x}+c\tilde{\gamma}^{i}=G^{\prime}(y)\tilde{\gamma}^{i}-\displaystyle\frac{1}{\mu_{i}}\sigma^{-2}\tilde{\psi}\chi_{\omega_{i}}& \mbox{ in }Q,\\[3mm]
\tilde{\psi}(0,t)=\tilde{\psi}(1,t)=0,~\tilde{\gamma}^{i}(0,t)=\tilde{\gamma}^{i}(1,t)=0& \mbox{ in }(0,T),\\[3mm]
\tilde{\psi}(x,T)=\tilde{\psi}_{T},~\tilde{\gamma}^{i}(x,0)=0& \mbox{ in }(0,1).
\end{array}\right.
\end{equation}
By the duality technique, it is easy to check that the null controllability of~(\ref{8.61})~with respect to the solution~$\tilde{z}$~can be reduced to prove the following observability inequality for (\ref{8.65}).

\begin{proposition}\label{p4}
~Assume that~$O_{d}\cap\omega\neq\emptyset,~x_{0}\in\omega$~and Condition~3~holds. If ~$\mu_{i}(i=1,2)$ are sufficiently large, then there exists a constant~$C>0$~such that for any~$\tilde{\psi}_{T}\in L^{2}(0,1)$~and~$y\in L^{2}(Q)$, the solution ~$(\tilde{\psi},\tilde{\gamma}^{i})$~of system~(\ref{8.65}) satisfies
\begin{eqnarray}\label{8.66}
\displaystyle\int_{0}^{1}|\tilde{\psi}(x,0)|^{2}dx+\displaystyle\sum_{i=1}^{2}\int_{Q}|\tilde{\gamma}^{i}|^{2}dxdt\leq C\int_{\omega\times(0,T)}|\tilde{\psi}|^{2}dxdt.
\end{eqnarray}
\end{proposition}

Since~(\ref{8.65})~is a linear system, similar to Theorem~\ref{ttt},~Proposition~\ref{p4}~can be obtained directly. Moreover, $\tilde{f}$~satisfies
\begin{eqnarray}\label{8.67}
J(\tilde{f})=\displaystyle\min_{f}J(f),~~\forall f\in L^{2}(\omega\times(0,T)).
\end{eqnarray}
Similar to~(\ref{8.20}), for $\mu_{i}(i=1,2)$ sufficiently large, it holds that 
\begin{eqnarray}\label{8.68}
\|\tilde{f}\|_{L^{2}(\omega\times(0,T))}\leq C\|\tilde{\psi}_{T}\|_{L^{2}(0,1)},~\forall\ y\in L^{2}(Q).
\end{eqnarray}
Therefore, we obtain that for any~$y\in L^{2}(Q)$,there exists  a control~$\tilde{f}\in L^{2}(\omega\times(0,T))$ such that the associated solution to system(\ref{8.61}) satisfies
\begin{eqnarray}\label{8.69}
\tilde{z}(\cdot,T;\tilde{f})=0~in~(0,1), \mbox{and}~(\ref{8.67})~\mbox{holds}.
\end{eqnarray}

Let's define a mapping~$\Pi: L^{2}(Q)\rightarrow L^{2}(Q)$,~with~$\Pi(y)=\tilde{z}_{y},~\forall y\in L^{2}(Q)$, where~$\tilde{z}_y$~is the solution of  system~(\ref{8.61})~corresponding to the control~$\tilde{f}$, that is,~$\tilde{z}_y$ satisfies~(\ref{8.69}). Obviously,~$\Pi$~is well defined. Combining~(\ref{8.64})~with~(\ref{8.68}), we get that $\tilde{z}$ and $\tilde{z}_{t}$  are uniformly bounded in $L^{2}(0,T;H_{a}^{1}(0,1))$ and $L^{2}(0,T;H^{*}_{a}(0,1))$, respectively. Therefore,  using Aubin-Lions~compactness theorem,   the mapping~$\Pi$ is compact  from~$L^{2}(Q)$~into itself. Also, $\Pi$ is continuous. Thus, by the~Schauder~fixed point theorem,  there is a fixed point of $\Pi$. This implies that, 
 we can find a control ~$\tilde{f}\in L^{2}(\omega\times(0,T))$~such that the solution of~(\ref{8.59})  satisfies (\ref{8.60}) and (\ref{8.888}). This ends the proof of Theorem \ref{T1} for the semilinear case. \endpf

\end{document}